\documentclass[letterpaper,11pt,reqno]{amsart}

\makeatletter
\usepackage{amssymb}
\usepackage{latexsym}
\usepackage{amsbsy}
\usepackage{amsfonts}
\usepackage{hyperref}
\usepackage[T1]{fontenc}
\usepackage{mathtools}   
\usepackage{cases}

\def\marginpar#1{\ignorespaces}

\newtheorem{theorem}{Theorem}[section]
\newtheorem{lemma}[theorem]{Lemma}
\newtheorem{proposition}[theorem]{Proposition}
\newtheorem{corollary}[theorem]{Corollary}

\newtheorem{definition}[theorem]{Definition}
\def\te{\rightarrow}
\def\xpsi{{\psi_0}}
\def\LK{{L\'evy-Khintchine }}
\def\psizero{{\psi_\alpha}}
\def\psione{{\tilde{\psi}_\alpha}}
\def\psioneone{{\tilde{\psi}_1}}
\def\lnu{{L}}
\def\sg{{c}}
\def\ed{ \stackrel{d}{=} }
\def\convd{ \stackrel{d}{\rightarrow} }
\def\sgn{ { \,{\rm sign} } }
\def\implies{\Rightarrow}
\def\hf{ \mbox{$\frac{1}{2}$} }
\def\pitwo{ \mbox{$\frac{\pi}{2}$} }
\def\twovpi{ \mbox{$\frac{2}{\pi}$} }

\def\reals { \mathbb{R} }

\def\mm{ \kappa }

\numberwithin{equation}{section}

\newdimen\AAdi%
\newbox\AAbo%
\def\AAk#1#2{\setbox\AAbo=\hbox{#2}\AAdi=\wd\AAbo\kern#1\AAdi{}}%

\def\eqref#1{(\ref{#1})}
\def\eqlabel#1{\def\@currentlabel{#1}}

\def\formula#1{\def\@tempa{#1}\let\@tempb\theequation\def\theequation{%
\hbox{#1}}\def\@currentlabel{(\theequation)}$$}
\def\endformula{\leqno\hbox{(\@tempa)}$$\@ignoretrue\let\theequation\@tempb}

\def\given{\hskip5\p@\relax\vrule\@width.4\p@\hskip5\p@\relax}

\newcommand{\open}[1]{%
\par\normalfont\topsep6\p@\@plus6\p@\trivlist\item[\hskip\labelsep\itshape#1%
\@addpunct{.}]\ignorespaces}

\DeclareRobustCommand{\close}[1]{%
  \ifmmode 
  \else \leavevmode\unskip\penalty9999 \hbox{}\nobreak\hfill
  \fi
  \quad\hbox{$#1$}}

\newlength{\toskip}

\makeatother

\begin{document}

\title [Stable distributions]{A direct approach to the stable distributions}

\author[E. J. G. Pitman] {{E. J. G.} Pitman}
\address{E. J. G. Pitman: 29 October 1897 $-$ 21 July 1993}
\author[Jim Pitman] {{Jim} Pitman}
\address{Jim Pitman: Statistics department, University of California, Berkeley} 
\email{pitman@stat.berkeley.edu}

\date{\today} 

\begin{abstract}
The explicit form for the characteristic function of a stable distribution on the line is derived analytically by solving 
the associated functional equation and applying theory of regular variation, without appeal to the general L\'evy-Khintchine 
integral representation of infinitely divisible distributions. 

\end{abstract}
\maketitle
\textit{Key words:} stable distribution, characteristic function, regular variation, infinitely divisible distribution, L\'evy-Khintchine formula, Eulerian integrals

\textit{AMS 2010 Mathematics Subject Classification: 60E07}
\section{Introduction}

This article is an expansion by Jim Pitman of a typescript with the same title which he received from his father E. J. G. Pitman in November 1983. 
The aim of the typescript was to derive the explicit form for the characteristic function of a stable distribution on the line, by solving 
the associated functional equation (\eqref{11} below)
and applying theory of regular variation, without appeal to the general L\'evy-Khintchine 
representation of infinitely divisible distributions, but with the necessary restrictions on the constants.
In a letter dated the 6th of November, 1983, my father wrote about the typescript:
\begin{quote}
I gave a rough version of the first four pages as long ago as when I was still at the University of Tasmania. I included an improved version in
my lectures at Chicago in 1969; but I was unable to get the restrictions on the constants. That I succeeded in doing only this year.
\end{quote}
To provide some context, this introduction and the following three sections offer a brief 
survey of the theory of stable distributions on the line and their L\'evy-Khintchine representations, based on the classical
theory of Eulerian integrals.
The remaining sections follow my father's 1983 typescript quite closely, with only minor changes of notation and arrangement of the material, and addition of
some references. 
The first four pages of the typescript, mentioned above,  became Section \ref{fnleq} 
of this article.
The argument then involves some extensions of my father's work \cite{pitman68} relating regular variation of  the tails of a distribution to regular variation of 
its characteristic function near the origin.
These results, which may be of some independent interest, are presented in Section  \ref{regvar}, followed by a proof of the main new result in 
Section \ref{regvarproof}.
The results on regular variation are finally applied in Section \ref{constantsec} to establish the restrictions on the  constants 
in the analytic form of the stable characteristic functions.

Let $X$ be a real-valued random variable with distribution  function $F(x) := P(X \le x)$ and
characteristic function 
\begin{align}
\phi(t):= E ( e^{i t X } ) = \int_{-\infty}^\infty e^{itx } d F(x) .
\end{align}
The distribution of $X$ is called {\em stable} if for each positive integer $n$,
with $X_1, \ldots, X_n$ independent random variables with the same distribution as $X$, 
there is the equality in distribution
\begin{align}
\label{s1}
X_1 + \cdots + X_n \ed a_n X + b_n
\end{align}
for some real constants $a_n$ and $b_n$ with $a_n > 0$.
See Section \ref{secstable} for some equivalent forms and variations of this definition.
We aim to prove:

\begin{theorem} 
\label{thmstable}
{\rm (L\'evy-Khintchine  representation of stable distributions
\cite{MR1504839},
\cite{khintchine-levy}
)}
\\
A function $\phi$ is the characterisic function of a stable probability distribution on the real line  if and only if 
$\phi(t) = \exp[\psi(t)]$  with 
\begin{align}
\label{lkstable}
\psi(t) = \psi_{\alpha,\sg,\beta,\mu}(t):= i \mu t -  |\sg t| ^\alpha ( 1 - i \beta \sgn(t) \omega_\alpha(t) \,)
\end{align}
where
\begin{subequations}
    \begin{numcases}
{
\omega_\alpha(t)  = 
}
\tan \, \pitwo \alpha  \mbox{ if } \alpha \ne 1 \label{omeganot1} \\
-\frac{2 }{\pi} \log |t|     \quad \mbox{ if } \alpha = 1  \label{omega1}
    \end{numcases}
\end{subequations}
\begin{align}
\label{constraints}
\mbox{ for some real }  \alpha, \sg, \beta, \mu \mbox{  with } 0 < \alpha \le 2, \sg \ge 0  \mbox{ and } |\beta| \le 1.
\end{align}
\end{theorem}
Here by definition, $\sgn(t) = t/|t| = 1( t >0 ) - 1 ( t < 0)$ has value $1$ or $0$ or $-1$ according to whether $t >0$ or $t = 0$ or $t <0$. 
This formulation of Theorem \ref{thmstable} follows the treatment by 
Gnedenko and Kolmogorov \cite[\S 34]{MR0062975}, 
with corrections indicated by Hall \cite{hall81}.  Following Samorodnitsky and Taqqu \cite[Definition 1.1.6]{MR1280932}, 
we encapsulate Theorem \ref{thmstable} as follows:
\begin{definition} 
\label{salparams}
\em{
The probability distribution whose characteristic function has logarithm $\psi_{\alpha,\sg,\beta,\mu}(t)$ is denoted 
$S_\alpha(\sg , \beta, \mu)$.
The notation $X \sim S_\alpha(\sg , \beta, \mu)$ signifies that a random variable $X$ has distribution $S_\alpha(\sg, \beta, \mu)$ for some $(\alpha, \sg, \beta, \mu)$ subject to the
constraints \eqref{constraints}. 
Then $X$, or the distribution of $X$, is called {\em stable with index $\alpha$}, or {\em $\alpha$-stable}.
Since $\omega_2(t) \equiv 0$, the value of $\psi_{2,\sg,\beta,\mu}(t)$ is unaffected by $\beta$. 
We adopt the convention that $\beta = 0$ for $\alpha = 2$. So $S_2(\sg, 0, \mu)$ is the normal distribution with mean $\mu$ and variance $2 \sg^2$.
}
\end{definition}
Hall's article \cite{hall81} reviewed and corrected a number of derivations of 
Theorem \ref{thmstable} 
in the literature up to 1980.
All these derivations 
follow the classical approach via the more general representation of distributions of $X$ that are {\em infinitely divisible}, meaning
that for every positive integer $n$
\begin{align}
\label{infdiv}
X \ed X_{n,1} + \cdots + X_{n,n}
\end{align}
for some sequence of $n$ independent and identically distributed random variables $X_{n,1}, \ldots X_{n,n}$. 
If the distribution of $X$ is stable as in \eqref{s1}, then \eqref{infdiv} holds with $X_{n,i} = (X_i - b_n/n)/a_n$.
Theorem \ref{thmstable} is then seen to be the specialization to the stable case of 

\begin{theorem} 
\label{thminfdiv}
{\rm (L\'evy-Khintchine  representation of infinitely divisible distributions
\cite{MR1504839},
\cite{khintchine-levy}
)
}
\label{infdivthm}
\\
A function $\phi$ is the characterisic function of an infinitely divisible probability distribution on the line  
if and only if $\phi(t) = \exp[ \psi(t) ]$ with 
\begin{align}
\label{lkinfdiv}
\psi(t) =  i b t - \hf \sigma^2 t^2 + \int_{-\infty}^\infty ( e^{i t x } - 1 - i t \tau(x) ) \lnu(dx)
\end{align}
where 
\begin{itemize}
\item $b$ is real
\item $\sigma \ge 0$, 
\item $\tau(x)$ is a bounded measurable {\em truncation function} such that $\tau(x)/x \te 1$ as $x \te 0$,
\item $\lnu$ is a {\em L\'evy measure} on the line, meaning that $\lnu(\{0\}) = 0$ and $\int_{-\infty}^\infty  (x^2 \wedge 1 )  \lnu(dx) < \infty$. 
\end{itemize}
\end{theorem}
Section \ref{secinfdiv} reviews details of the L\'evy-Khintchine representation \eqref{lkinfdiv} in the particular case of stable laws, and
the derivation of the ``if'' part of Theorem \ref{thmstable}, that the distribution $S_\alpha(\sg, \beta, \mu)$ is well defined, 
from the ``if'' part of Theorem \ref{infdivthm}, which is relatively easy to prove compared to its ``only if'' part. 
This derivation in Section \ref{secinfdiv} relies on some classical evaluations of Eulerian integrals, which are reviewed in  Section \ref{eulerian}.
The ``only if'' part of Theorem \ref{thmstable}, that every stable distribution is of the form $S_\alpha(\sg, \beta, \mu)$ for some $(\alpha, \sg, \beta, \mu)$ 
subject to the constraints \eqref{constraints},
takes much more work. 
For many decades, there seemed to be no alternative to first establishing the ``only if'' part of the general L\'evy-Khintchine representation, which involves 
a fair amount of analytic machinery, including compactness arguments, the Helly-Bray theorem and L\'evy's continuity theorem for characteristic functions,
then specializing to obtain the form of the L\'evy measure of an infinitely divisible distribution that is stable.

The work of Geluk and de Haan \cite{geluk-haan} in 2000  seems to be the first published alternative to this classical approach.
Their approach is similar to that adopted here, using Fourier transforms and theory of regular variation.
Their treatment is also more complete. They start from the characterization of stable laws of $X$ as limits in distribution 
of centered, scaled sums $X_1 + \cdots + X_n$ of independent and identically distributed $X_i$, say
\begin{align}
\label{limd}
(X_1 + \cdots + X_n - \mu_n)/c_n \convd X ,
\end{align}
where all the $X_i$ have the same distribution as $X_1$, but the limit distribution of $X$ may differ from the
distribution of $X_1$. For instance, in the most familiar form of the central limit theorem, 
starting from any distribution of $X_1$ with mean $\mu$ and variance $\sigma^2$,
 leads with $\mu_n = n \mu$ and $c_n = \sqrt{n} \sigma$ to $X$ with standard normal distribution.
Geluk and de Haan go on to describe the {\em domain of attraction} of each possible stable law of $X$, that
is the set of all possible distributions of $X_1$ such that \eqref{limd} holds for independent and identically distributed $X_i$, 
along with the associated constants $\mu_n$ and $c_n$. 
In this approach, the L\'evy-Khintchine representation of stable laws emerges as a byproduct of a comprehensive analysis of their domains of attraction.
See also Feller \cite[Section XVII.5]{fellerII},
Mijneer \cite{MR0370783}, 
and the references of \cite{geluk-haan} 
for other approaches to domains of attraction.

The present approach to the ``only if'' part of Theorem \ref{thmstable} is both more direct and less ambitious 
than the approach of Geluk and de Haan. It involves only an analysis of the functional equation satisfied by the characteristic function of stable distribution, \eqref{11} below, without consideration of domains of attraction.
See also the text of Ramachandran and Lau \cite[Chapter 3]{MR1132671} for a similar analysis of the functional equation, 
with references to earlier work.
For general background on regular variation and its applications to probability theory, see the texts of Feller \cite{fellerII}, Bingham, Goldie and Teugels \cite{bgt}, 
and Geluk and de Haan \cite{MR906871}, and the historical article of Seneta \cite{MR1967900}.
See also Kagan, Linnik and Rao \cite{MR0346969} regarding characterization problems in probability and statistics,
Kuczma \cite{MR2467621} and 
Acz{\'e}l and Dhombres  \cite{MR1004465}
for general background on functional equations, 
and 
\cite{MR1491578}, 
\cite{MR2213594} 
for some more recent work on characterizations of multivariate stable laws through functional equations satisfied by their characteristic functions.

\section{Stable distributions}
\label{secstable}
The definition of stability \eqref{s1} can be written in terms of distribution functions as
\begin{align}
\label{s2}
F^{n*}(x) = F( (x - b_n)/a_n)
\quad \mbox{ for all real } x ,
\end{align}
where $F^{n*}$ is the $n$-fold convolution of $F$ with itself. 
Equivalently, 
in terms of characteristic functions,
\begin{align}
\label{11}
\phi(t)^n = \phi( a_n t ) \exp ( i b_n t )
\quad \mbox{ for all real } t .
\end{align}
Starting in Section \ref{fnleq}, we take the functional equation \eqref{11} satisfied by $\phi$ as the defining
property of stability, then proceed to characterize all solutions of this functional equation that are characteristic functions of some probability distribution on the line.
But a number of variations and generalizations of the stability concept also deserve mention. And there is no doubt that in the applications of stable laws,
their 
representation as infinitely divisible laws, especially their description as stochastic integrals with respect to Poisson processes governed
by the stable L\'evy measure \cite[Chapter 3]{MR1280932}, is of utmost importance.

First of all, there is the characterization of stable distributions as limits in distribution of centered and scaled sums of independent and
identically distributed random variables, as in \eqref{limd}. Then there is the condition that for $X_1$ and $X_2$ two independent copies
of $X$, for every pair of positive reals $a$ and $b$, there exists a positive real $c$ and a real $d$ such that
\begin{align}
\label{abcd}
aX_1 + b X_2 \ed cX + d .
\end{align}
The equivalence of these three variations of the definition of stability is elementary, and discussed in many texts. See for instance 
Feller \cite[VI.1]{fellerII}.
See also Samorodnitsky and Taqqu \cite{MR1280932}, 
for a comprehensive account of the theory of stable distributions and its applications to the construction of stable non-Gaussian random processes.

The following byproduct of Theorem \ref{thmstable}  is established by elementary arguments in Section \ref{fnleq}:
\begin{corollary}
\label{scaling}
The only possible form of the constants $a_n$ in \eqref{s1} is $a_n = n^{1/\alpha}$ for some $0 < \alpha \le 2$.
Then in \eqref{abcd} 
\begin{align}
a^\alpha + b^\alpha = c^\alpha  \qquad \qquad ( 0 < \alpha \le 2 ).
\end{align}
\end{corollary}
The distribution of $X$ is called {\em symmetric} if $X \ed -X$ and 
{\em strictly stable} if any one of the following equivalent conditions obtains:
\begin{itemize}
\item \eqref{abcd} holds with $d = 0$;
\item \eqref{s1} holds with $b_n = 0$;
\item \eqref{limd} holds with $d_n = 0$.
\end{itemize}

By inspection of the characteristic function \eqref{lkstable}, and appeal to the uniqueness theorem for characteristic functions, 
$X$ is 
\begin{itemize}
\item 
symmetric $\quad \,\, \Leftrightarrow$ either $(\alpha = 2$  and  $\mu =0)$ or $(0 < \alpha < 2$ and $\mu = 0$ and $\beta = 0)$;
\item 
strictly stable $\Leftrightarrow$ either $(0 < \alpha \le 2$ and $\mu = 0)$ or $(\alpha = 1$ and $\beta = 0)$.
\end{itemize}
Note the subtle difference between the cases $\alpha \ne 1$ and $\alpha = 1$:
\begin{itemize}
\item if $\alpha \ne 1$, then for every $\alpha$-stable $X$, no matter what the value of $\beta$, there is a constant $\mu$ such the shift $X-\mu$ is strictly stable;
\item but if $\alpha = 1$, then distribution of $X$ is strictly stable for $\beta = 0$,  but for $\beta \ne 0$ no shift of $X$ is strictly stable.
\end{itemize}
This difference is indicative of a discontinuity  in the parameterisation of $S_\alpha(\sg , \beta, \mu)$
at $\alpha = 1, \beta \ne 0$, discussed further in the next section.  
The same difference is 
apparent in the following straightforward corollary of \eqref{lkstable}:
\begin{corollary}
\label{sumiid}
If $X \sim S_\alpha(\sg , \beta, \mu)$ and $S_n:= X_1 +  \cdots + X_n$ is the sum of $n$ independent copies $X_i$ of $X$, then
$S_n \sim S_\alpha(n^{1/\alpha} \sg , \beta, n \mu)$. Equivalently, 
\begin{subequations}
    \begin{numcases}
{
S_n  \ed 
}
n^{1/\alpha} X + \mu(n - n^{1/\alpha})  \mbox{ if } \alpha \ne 1 \label{n1}\\
n X + \frac{2}{\pi} c \beta n \log n             \quad \mbox{ if } \alpha = 1. \label{n2}
    \end{numcases}
\end{subequations}
\end{corollary}

It is a well known consequence of Theorem \ref{thmstable}, by Fourier inversion, that
each non-degenerate stable law has a bounded and continuous density. For instance, in the strictly stable case, up to a change of location and scale, the
density of $X$ can be written as follows \cite[Theorem V.7.13]{MR2011862}: 
\begin{align}
\label{ftheta}
f_\theta(x) = \frac{1}{\pi} \int_0^\infty \exp( - u ^\alpha ) \cos ( x u + \theta u^\alpha ) \, du \quad \mbox{ for all real } x 
\end{align}
where $\theta$  is a real parameter subject to $|\theta| \le |\tan \pitwo \alpha |$.
Unfortunately, there seem to be few more explicit formulas for $f_\theta$ except in familiar special cases: the Gaussian case $\alpha = 2$, the Cauchy case
$\alpha = 1$, and the one-sided L\'evy case with $\alpha = 1/2$ which arises from first passage times of Brownian motion.
See Uchaikin and Zolotarev \cite{MR1745764} for a deeper study of stable densities.
As remarked in \cite[\S (7.26)]{MR2011862}, Theorem \ref{thmstable} implies that for an integrable function of the form \eqref{ftheta} for some real $\theta$,
\begin{align}
\label{ftheta1}
f_\theta(x) \ge 0 \mbox{ for all real } x \Leftrightarrow  |\theta| \le |\tan \pitwo \alpha |.
\end{align}
This seems to be very hard to prove directly: so considerable work must be done to establish the limitation $|\beta| \le 1$ in the parameterization \eqref{lkstable} of stable laws.

The terminology used here has evolved over time. L\'evy used the terms {\em stable} and {\em quasi-stable}, instead of {\em strictly stable} and {\em stable}.
L\'evy also introduced the term {\em semi-stable} for distributions of $X$ such that \eqref{s1} holds for $n=2$, or equivalently \eqref{abcd} holds for $a=b$.
L\'evy characterized all distributions with this property, which is weaker than stability, but which can nonetheless be approached via the theory of infinitely divisible distributions. 
Another closely related concept is that of a {\em self-decomposable} distribution, or distribution of {\em class $L$}.
See Sato \cite{MR1739520} 
for a comprehensive treatment of all these notions, their generalizations to higher dimensions, and their applications to the
theory of stochastic processes with independent increments.
Steutel and van Harn \cite[Chapter V]{MR2011862} 
offers a more elementary treatment of self-decomposable and stable distributions on the line, including a proof of 
Theorem \ref{thmstable}
using a representation of self-decomposable distributions.
See also Bertoin \cite{MR1406564} regarding L\'evy processes.
Another research monograph with extensive treatment of stable distributions and their applications is 
Uchaikin and Zolotarev \cite{MR1745764}. 

\section{Eulerian integrals}
\label{eulerian}

The Fourier transforms of probability distributions on the line, especially those of stable distributions, are closely related to 
various Eulerian integrals. These definite integrals were first evaluated informally by Euler in the 18th century, then studied with increasing degrees of rigour in 
the 19th century by Legendre, Cauchy, Dirichlet, Saalsch\"utz and others. 
See Whittaker and Watson \cite[XII]{whittwat} for an account of Eulerian integrals from the viewpoint of classical analysis, with references to original sources.

Starting from Legendre's definition of the gamma function for positive real $r$
\begin{align}
\label{gamdef}
\Gamma(r) := \int_0^\infty t^{r-1} e^{-t} dt  \qquad  ( r > 0 )
\end{align}
there is the basic recursion
\begin{align}
\label{gamrec}
\Gamma(r + 1) = r \Gamma(r)
\end{align}
obtained by integration by parts for $r > 0$. The definition of $\Gamma(r)$ can then be extended to all real $r$
except $r = 0, -1, -2, \ldots$, either by repeated use of the recursion \eqref{gamrec}, or by Euler's reflection formula
\begin{align}
\label{reflection}
\Gamma(r) \Gamma(1 - r) = \frac{ \pi }{ \sin \pi r }.
\end{align}
While the definition of $\Gamma(r)$ for complex $r$ by analytic continuation has added profoundly to  theory of the gamma and related functions \cite{whittwat}, 
the focus here is on the classical integral representations of $\Gamma(r)$ for real $r$, especially $r = - \alpha$ with $0 < \alpha < 1$ or $1 < \alpha <2$,
which cases are of special importance for the evaluation of stable characteristic functions.
The following theorem presents these classical integral representations of $\Gamma(r)$ for $r < 1$.

\begin{theorem}
\label{cauchythm}
{\rm (Euler-Cauchy-Saalsch\"utz)}
For each real number $r < 1$ except $0, -1, -2 , \ldots$,
and each non-zero complex number $z$ with $\Re(z ) \ge 0$, 
\begin{align}
\label{cauchysaal}
\int_0^\infty x^{r-1} \left(e^{- z x } - \sum_{0 \le i < - r } \frac{ ( - z x )^i }{ i! } \right)  dx = \Gamma(r) z^{-r}  \quad 
\end{align}
\end{theorem}
Here the {\em compensating sum } $\sum := \sum_{0 \le i < - r }$ required to make the integral converge  is
\begin{itemize} 
\item for $ 0 < r <1$: an empty sum $\sum = 0$;  
\item for $r$ with $-(k + 1) < r < -k \le 0$ for some non-negative integer $k$: \\
$\sum = \sum_{0 \le i \le k}$ is the sum of the first $k + 1$ terms in the Maclaurin series of $e^{-zx}$.
\end{itemize}
If $\Re(z) >0$ the identity \eqref{cauchysaal} holds for all real $r$ except for $r = 0, -1, -2 , \ldots$, with an absolutely convergent integral on the left, and
$z^{-r}$ on the right defined by its principal value, that is
\begin{align}
\label{argdef}
z^{-r} := |z|^{-r} e^{- i r \arg(z) } = 
\frac{ e^{i r \arctan(\theta/\lambda)  }  }
{ (\lambda^2 + \theta^2)^{r/2} }
\mbox{ if } z = \lambda - i \theta \mbox{ for } \lambda >  0  , \theta \in \reals.
\end{align}
where $\arg(z) \in [-\pitwo, \pitwo]$ is the principal value of the argument of $z$ with $\Re z \ge 0$.  
So the simplest case of \eqref{cauchysaal} for $r >0$  with $\sum = 0$
for 
$z = \lambda - i \theta $ with $\lambda >  0 , \theta \in \reals$ reads
\begin{align}
\label{gamcf}
\int_0^\infty x^{r-1} e^{-\lambda x } e^{i \theta x } \, dx = 
\frac{ \Gamma(r)  }{ (\lambda - i \theta)^r } = 
\Gamma(r)  
\frac{ e^{i \arctan(\theta/\lambda) r } }
{(\lambda^2 + \theta^2 ) ^{r/2} }
\qquad( r > 0 , \lambda > 0, \theta \in \reals).
\end{align}
For fixed $r>0$ and $\lambda >0$, this function of $\theta$, multiplied by $\lambda^r/\Gamma(r)$, is the characteristic function of the gamma distribution 
with parameters $(r,\lambda)$ on the positive half line.  
As indicated by Feller in \cite[p. 502]{fellerII}, this formula follows easily from the series expansion of $e^{i \theta x }$ in powers of $x$
by integrating term by term.

For $r < 1$ and $\Re(z) \ge 0$ the integral in \eqref{cauchysaal} is an improper integral $\int_0^\infty:= \lim_{T \to \infty} \int_0^T$, 
which may or may not be absolutely convergent, depending on the choice of $r$ and $z$. 
Let $I(z,r)$ denote the value of this integral. The key observation is that for $r <0$, provided $r$ is not a negative integer, and $z \ne 0$ has $\Re(z) \ge 0$,
integration by parts gives
\begin{align}
\label{cauchysaalrec}
I(r,z) r/z = I(r+1,z) ,
\end{align}
hence easily $I(r,z) = \Gamma(r) z^{-r}$ by reduction to the case $r >0$, using the gamma recursion \eqref{gamrec} to define $\Gamma(r)$ for $r < 0$.
In particular, for $z = - i\theta$, with $\theta$ real and non-zero, on the right side of \eqref{cauchysaal} we see
\begin{align}
\label{cauchysaal1}
 z^{-r} = ( - i \theta )^{-r} =  |\theta|^{-r} e^ { i \sgn(\theta) \pitwo r }  = |\theta|^{-r} \left( \cos \pitwo r  + i \sgn(\theta) \sin \pitwo r \right) .
\end{align}
So \eqref{cauchysaal} gives for all real $\theta \ne 0$ and non-integer $r < 1$ 
\begin{align}
\label{cauchysaali}
\int_0^{\infty} x^{r-1} ( e^{i \theta x }   - \Sigma_{r, \theta x } ) dx   = \Gamma(r) |\theta|^{-r} \left( \cos \pitwo r  + i \sgn(\theta) \sin \pitwo r \right)
\end{align}
where the compensating sum $\Sigma_{r, \theta x }$ is
\begin{align*}
\Sigma_{r, \theta x } &= 0   \qquad \qquad \mbox{ if } 0 < r < 1  \\
&= 1   \qquad \qquad \mbox{ if }  -1 < r < 0   \\
&= 1 + i \theta x \quad \,\,\mbox{ if }  -2 < r < -1 
\end{align*}
and so on.
For general non-integer $r < 1$, the real and imaginary parts of \eqref{cauchysaali} can be unpacked as
\begin{align}
\label{rlpart}
\int_0^\infty x^{r-1} \left( \cos \theta x \,\, - \sum_{0 \le k < - r/2} \frac{ (-1)^k ( \theta x )^{2k}  } {(2k)!} \right)  dx = \Gamma(r) |\theta|^{-r} \cos \pitwo r
\end{align}
and 
\begin{align}
\label{impart}
\int_0^\infty x^{r-1} \left( \sin \theta x \,\, - \sum_{0 \le k < -(r+1)/2} \frac{ (-1)^k (\theta x )^{2k + 1} } {(2k +1)!} \right) dx  = \Gamma(r) |\theta|^{-r} \sgn(\theta) \sin \pitwo r  .
\end{align}
The following corollary of Theorem \ref{cauchythm} provides some useful variants of these formulas:
\begin{corollary}
\label{sincos}
For each positive real $\mm$ that is not an odd integer,  and all real $\theta \ne 0$
\begin{align}
\label{rlpart1}
\int_0^\infty x^{- \mm}  \left( \cos \theta x \,\, - \sum_{0 \le j < (\mm - 1)/2} \frac{ (-1)^j ( \theta x )^{2j} } {(2j)!} \right)  dx = C( \mm) | \theta|^{\mm - 1 }
\end{align}
where 
\begin{align}
\label{cm}
C(\mm) =  {  \Gamma( 1 - \mm) \sin( \pitwo \mm )} = 
\frac{ \pitwo }{  \Gamma(\mm) \cos( \pitwo \mm )},
\end{align}
and for each positive real $\mm$ that is not an even integer, and all real $\theta \ne 0$
\begin{align}
\label{impart1}
\int_0^\infty x^{- \mm}  \left( \sin \theta x \,\, - \sum_{0 \le j < (\mm - 2)/2} \frac{ (-1)^j ( \theta x )^{2j + 1} } {(2j +1)!} \right)  dx = S( \mm) \sgn(\theta) | \theta|^{\mm - 1 }
\end{align}
where
\begin{align}
\label{sm}
S(\mm) = 
{  \Gamma( 1 - \mm) \cos ( \pitwo \mm )} =
\frac{ \pitwo }{ \Gamma(\mm) \sin( \pitwo \mm )}
\end{align}
so that
\begin{align}
\label{smcmrat}
\frac{C(\mm)} {S(\mm)} =  \tan (\pitwo \mm) .
\end{align}
\end{corollary}
These formulas, at first for positive non-integer $\mm$ with the first expressions for $C(\mm)$ and $S(\mm)$ involving $\Gamma(1-\mm)$, 
and the consequent ratio \eqref{smcmrat},
 are read from \eqref{rlpart} and \eqref{impart} for $r = 1 - \mm$, using 
$$
\cos \pitwo (1- \mm) = \sin \pitwo \mm \quad \mbox{and} \quad \sin\pitwo (1 - \mm) = \cos\pitwo \mm. 
$$
The alternative expressions for $C(\mm)$ and $S(\mm)$ involving $\Gamma(\mm)$ are then read from Euler's reflection formula \eqref{reflection} and the sine 
duplication formula $\sin (2 \theta) = 2 \sin \theta \cos \theta$ for $\theta = \mm \pi/2$.
The formulas for the cosine integral \eqref{rlpart} with $\mm$ an even integer, with $C(\mm)$ evaluated by continuity using the expression with $\Gamma(\mm)$,
and  for the sine integral \eqref{impart1} for $\mm$ an odd integer, with similar evaluation of $S(\mm)$, are obtained by analytic continuation.

The case $\mm = 1$ of the sine integral \eqref{impart1}, with $S(1) = \pi/2$,  is of particular importance. This is the {\em Dirichlet integral} \cite[p. 315]{courantv2} 
\begin{align}
\label{dirichlet}
\int_0^{\infty} \frac{ \sin \theta x }{x} dx = \sgn(\theta) \frac{\pi}{2}  
\end{align}
which is the basis of the Fourier inversion formulas for the characteristic function of a probability distribution, which are
recalled in Section \ref{regvar}. 
See also von Bahr \cite{MR0179827} for an application of \eqref{rlpart1} to give formulas for fractional moments in terms of 
the characteristic function, and Pitman \cite{pitman68} for further applications to the  analysis of characteristic functions.
\section{Infinitely divisible laws}
\label{secinfdiv}

To establish the 
the existence of the distribution denoted $S_\alpha(\sg,\beta,\mu)$ in Definition \ref{salparams},
that is the ``if'' part of Theorem \ref{thmstable}, 
it must be shown that each of the functions presented in \eqref{lkstable} is the characteristic function of a probability distribution on the line. 
For that done, it follows  easily that the distribution $S_\alpha(\sg,\beta,\mu)$ is stable with index $\alpha$,
as indicated already in Corollary \ref{sumiid}.

For $\alpha = 2$, the distribution $S_2(\sg,\beta,\mu)$ is just normal with mean $\mu$ and variance $2 \sg^2$.
Stability of the normal distribution follows by evaluation of its characteristic function, or by the convolution formula for normal densities.

For $0 < \alpha < 2$, following L\'evy, 
a random variable $X$ with $S_\alpha(\sg,\beta,\mu)$ distribution can be exhibited as a limit of centered compound Poisson variables.  
More generally, it is straightforward to show that for every triple $(b,\sigma^2,L)$ as in Theorem \ref{infdivthm}, 
\begin{itemize}
\item the function $\phi$ defined by \eqref{lkinfdiv} is the characteristic function of an infinitely divisible probability distribution on the line, and 
\item a random variable $X$ with this distribution can be constructed as the limit in distribution of a suitably centered sequence of compound Poisson variables.
\end{itemize}
This argument appears in many places, see  e.g. Feller \cite[XVII.2]{fellerII},  Kallenberg \cite[Corollary 15.8]{kallenberg02} or
Durrett \cite[(7.7)]{MR2722836}, 
so will not be repeated in detail here. 

To provide a complete account of the \LK representation of stable distributions, with correct values of the constants, it is necessary to specify the choice 
of truncation function $\tau(x)$ appearing in the general \LK representation \eqref{lkinfdiv}.
Some common choices are
\begin{align}
\tau(x) &:=  x/(1+x)^2  &\mbox{ L\'evy \cite{levy54}, Gnedenko-Kolmogorov \cite{MR0062975}} \\
\tau(x) &:=  \sin x   &\mbox{ Feller [\S XVII]\cite{fellerII} }\\
\tau(x) &:= x 1(|x| \le 1 ) &\mbox{ Kallenberg \cite[Theorem 15.9]{kallenberg02} }
\end{align}
For a particular infinitely divisible characteristic function whose logarithm $\psi$ has the L\'evy-Khintchine representation \eqref{lkinfdiv}, whatever the choice of $\tau$, the integral in 
the \LK representation is absolutely convergent. 
The difference between the integrals associated with two different truncation functions, say $\tau_1$ and $\tau_2$, is just $ict$ for the integral
\begin{align}
c = \int_{- \infty}^\infty ( \tau_1(x) - \tau_2(x) ) L(dx) 
\end{align}
which is absolutely convergent by the properties of $\tau_1, \tau_2$ and $L$ listed in Theorem \ref{thminfdiv}. So replacement of $\tau_1$ by $\tau_2$ just increments $b$ by $c$.
As Hall \cite{hall81} observes, the passage from \eqref{lkinfdiv} to \eqref{lkstable} for the stable L\'evy measure (displayed in \eqref{matchup} below) requires some care to avoid errors in the sign of $\beta$.
This argument is made unnecessarily complicated in many sources by the common choice of truncation functions other than the choice $\tau(x) = \sin x$
advocated by Feller. 
Use of $\tau(x) = \sin x$ greatly simplifies computation of the L\'evy-Khintchine integrals for the $\alpha$-stable L\'evy measures, as indicated in the
following corollary, which combines and simplifies presentations of Feller \cite[\S XVII]{fellerII}
and Zolotarev \cite[(M) on p. 11]{MR854867}.

\begin{corollary}
\label{matchmore}
Suppose $0 < \alpha < 2$.
For $X$ with the $\alpha$-stable distribution $S_\alpha(\sg , \beta, \mu)$, let
\begin{subequations}
\label{blabel}
    \begin{numcases}
{
b =
}
\mu + \beta \sg^\alpha \tan \pitwo \alpha \mbox{ if } \alpha \ne 1 \label{not1} \\
\mu               \quad \quad \quad \mbox{ if } \alpha = 1.  \label{eq1}
    \end{numcases}
\end{subequations}
Then from {\em \eqref{lkstable}} the logarithm of the characteristic function  of $X-b$ is
\begin{subequations}
\label{psialsgbeta}
    \begin{numcases}
{
\psi_{\alpha,\sg,\beta } (t) :=
}
-   \sg^\alpha \left( |t|^\alpha  - i \beta \tan ( \pitwo \alpha ) ( |t|^\alpha - 1 )  \right) \mbox{ if } \alpha \ne 1 \label{not11} \\
-  \sg  \left(  |t|  + i \twovpi t \log |t|   \right) \quad \quad \quad \mbox{ if } \alpha = 1 . \label{eq1}
    \end{numcases}
\end{subequations}
This function admits the \LK integral representation
\begin{align}
\label{lkcont}
\psi_{\alpha, \sg, \beta}(t) = \int_{-\infty}^\infty ( e^{i t x } - 1 - i t \sin x) \lnu_{\alpha,\sg,\beta}(dx) 
\end{align}
as in {\em \eqref{lkinfdiv}}, with truncation function $\tau(x) = \sin x$ and $b = \sigma^2= 0$, for the L\'evy measure
\begin{align}
\label{matchup}
\lnu_{\alpha, \sg, \beta}(dx) = \left[ \, (1-\beta)  1(x < 0 ) + (1 + \beta )1 (x >0 ) \,\right] \, \frac{ \sg^\alpha \, dx } { 2 K(\alpha) |x|^{1 + \alpha} } 
\end{align}
where 
\begin{align}
\label{kalalpha}
K(\alpha):= - \Gamma(-\alpha) \cos \pitwo \alpha = \frac{ \Gamma(2 - \alpha) }{\alpha (1 - \alpha) } \cos \pitwo \alpha 
\end{align}
is a continous and strictly positive function of $\alpha$ for $0 < \alpha < 2$, with  evaluation by
continuity at $\alpha = 1$:
 \begin{align}
\label{kal1}
K(1):= \pitwo .
\end{align}
\end{corollary}

As observed by Zolotarev \cite[(M) on p. 11]{MR854867}, 
the function
$\psi_{\alpha, \sg, \beta}(t)$ in \eqref{lkcont} and  \eqref{matchup} is a continuous function of $\alpha, \sg, \beta, t$ as $(\alpha,\sg, \beta)$ ranges over
the allowed parameter space 
$$
\mbox{$0 < \alpha < 2, \sg \ge 0, |\beta| \le 1$ and $-\infty < t < \infty$.}
$$

For $X \sim S_\alpha(\sg , \beta, \mu)$ and
$b = b_\alpha(\sg , \beta, \mu)$ as in \eqref{blabel}, the distribution of $X-b$ is a weakly continuous function of $(\alpha,\sg, \beta)$.
However, the shift $b$ is a spectacularly discontinuous function of $\alpha, \sg, \beta,\mu$ at $\alpha = 1, \beta \ne 0$:
as $\alpha$ increases to $1$ the factor of $\tan \pitwo \alpha $ explodes up to $+\infty$, while as $\alpha$ decreases to $1$ the same factor explodes down to $-\infty$, and 
the value of $b$ at $1$ is just $\mu$. Thus the parameterization of $S_\alpha(\sg , \beta, \mu)$ is extremely discontinuous at $\alpha = 1, \beta \ne 0$.
See Samorodnitsky and Taqqu \cite[p. 7,(1.1.10)] {MR1280932} 
for further discussion of this issue.
See also Geluk and de Haan \cite[Theorem 1]{geluk-haan} for a description of the domain of attraction of each stable distribution, using the 
distribution whose characteristic function has logarithm \eqref{psialsgbeta} for $\sg = 1$ as a convenient representative of the equivalence class of all stable 
distributions of a given type (meaning equivalent up to a change of location and scale).

Corollary \ref{matchmore} follows easily from the integral evaluations of the following lemma, as indicated in
Feller \cite[\S XVII.3 pp. 568-570]{fellerII} and Zolotarev \cite[(1.17) p.9]{MR854867}.

\begin{lemma}
\label{lemmalk}
For $0 < \alpha < 2$ and $K(\alpha)$ as in {\em \eqref{kalalpha}-\eqref{kal1}}, let
\begin{align}
\psizero (t):= - K(\alpha) |t|^{\alpha} ( 1 - i \sgn(t) \tan(\pitwo \alpha) ).
\end{align}
Then there are the integral representations
\begin{align}
\psizero (t) &=  \int_0^\infty( e^{itx} - 1 ) \frac{dx}{x^{\alpha + 1} } & \mbox{     if    } 0 < \alpha < 1 \label{case1}\\
&=  \int_0^\infty( e^{itx} - 1 - it x ) \frac{dx}{x^{\alpha + 1}} &                \mbox{     if    } 1 < \alpha < 2 \label{case2}
\end{align}
Moreover, the L\'evy-Khintchine characteristic exponent 
\begin{align}
\label{psionelk}
\psione(t):= \int_0^\infty( e^{itx} - 1 - i t \sin x ) \frac{dx}{x^{\alpha + 1}}  \qquad \mbox{ for } 0 < \alpha < 2
\end{align}
is jointly continuous in $(t,\alpha)$, and given by
\begin{align}
\psione(t) &= - K(\alpha) [ \, |t|^\alpha - i t \tan( \pitwo \alpha   ) ( |t|^{\alpha -1 } - 1 )\,  ] & \mbox{ for } \alpha \ne 1,  \mbox{ and } \label{psione}\\
\psioneone(t)   &= - \pitwo |t| - i t \log |t|                   				       & \mbox{ for } \alpha = 1. \quad  \quad \label{limpsione} 
\end{align}
\end{lemma}
\textbf{Proof:} 
As recognized by 
L\'evy \cite[\S 57, (35') and  p. 200, Deuxieme cas]{levy54} 
the integral representations \eqref{case1} and 
\eqref{case2} are the instances $r= - \alpha$ of the Eulerian integrals \eqref{cauchysaali}.
See also Feller \cite[\S XVII.3, p. 568]{fellerII} and 
Samorodnitsky and Taqqu \cite[Exercise 3.14, p. 170]{MR1280932}
for other derivations of these definite integrals. 
The equality of imaginary parts of \eqref{case1} for $t = 1$ gives
\begin{align}
\int_0^\infty \sin x \frac{dx}{x^{\alpha + 1} } = \Im \psizero(1) \quad (0 < \alpha < 1).
\end{align}
This yields \eqref{psionelk} for $0 < \alpha < 1$ with $\psione(t) = \psizero(t) - i t \, \Im \psizero(1)$ as in \eqref{psione}.
Similarly, \eqref{case2} gives
\begin{align}
\int_0^\infty (\sin x  - i x ) \frac{dx}{x^{\alpha + 1} } = \Im \psizero(1) \quad (1 < \alpha < 2),
\end{align}
hence \eqref{psionelk} for $1 < \alpha < 2$. It is easily shown using
\eqref{psionelk}
that $\psione(t)$ is jointly continuous in $(t,\alpha)$. So \eqref{limpsione} for $\alpha = 1$ 
is obtained in the limit, using \eqref{kal1},
$$
\lim_{\alpha \te 1}  (1-\alpha) \tan( \pitwo \alpha   )  = \frac{2}{\pi} \mbox{  and } \lim_{\alpha \te 1 } \frac{ t^{\alpha - 1} - 1 }{\alpha - 1 } = \log | t | .
$$
As a check, see also Zolotarev \cite[(1.17) p.9]{MR854867} 
where $Q_\alpha = \alpha \psione(t)$ and the formulas \eqref{psione} and \eqref{limpsione} appear in the middle of the following page.
$\square$

For the function $\psi_{\alpha, \sg, \beta}(t)$ defined by \eqref{psialsgbeta}, 
the integral representation \eqref{lkcont} follows immediately, by taking a linear combination of contributions from the positive and negative half-lines.
To complete the proof of existence of the $S_\alpha(\sg , \beta, \mu)$ distribution, it only remains to show that
$\psione(t)$ is the logarithm of a characteristic function for each $0 < \alpha < 2$. But it is well known and easily proved that this is the limit of 
the  logarithm of the characteristic function of a suitably centered sequence of compound Poisson variables which converges in probability, and the conclusion follows.
See Samorodnitsky and Taqqu \cite[Theorem 3.12.2]{MR1280932} for a comprehensive account of the Poisson representation of stable random variables, and
numerous other representations of  random variables $X \sim S_\alpha(\sg , \beta, \mu)$.

Following the classical approach of L\'evy and Khintchine, the  
``only if'' part of Theorem \ref{thmstable} can now be deduced from the ``only if'' part of the L\'evy-Khintchine representation 
infinitely divisible distributions (Theorem \ref{infdivthm}). 
See Gnedenko and Kolmogorov \cite[\S 34]{MR0062975}, 
other sources cited in Hall \cite{hall81},
or the presentation of Kallenberg \cite[Theorem 15.9]{kallenberg02} in the broader context of stable L\'evy processes, where
Theorem \ref{infdivthm} can be read from \cite[Corollary 15.8 and Theorem 15.12]{kallenberg02}.

\section{Analysis of the functional equation}
\label{fnleq}
Suppose throughout this section that $X$ is a random variable with a stable distribution. So the sum, $S_n$, of $n$ independent random variables,
each with the same distribution as $X$,
is distributed like $a_n X + b_n$, where $a_n$, $b_n$ are real, $a_n >0$.
If $\phi(t)$ is the characteristic function of $X$, the characteristic function of $S_n$ is
\begin{align}
E \left[ \exp ( i t [ a_n X + b_n] ) \right] = \phi(a_n t ) \exp ( i b_n t ).
\end{align}
Thus
\begin{align}
\label{1}
\phi(t)^n = \phi( a_n t ) \exp (i b_n t ).
\end{align}
Here and in the following, $m$ and $n$ will always denote positive integers.
To avoid trivialities, we will assume throughout that the distribution of $X$ is non-degenerate, meaning it is not concentrated at a single point. 
Equivalently, in terms of the characteristic function,
\begin{align}
\label{notdegen}
|\phi(t)| < 1 \mbox{  for some real } t. 
\end{align}
This section derives some first consequences of the functional equation \eqref{1}, including the following proposition.
See also L\'evy \cite[\S 95]{levy54}, Feller [Theorem VI.1.1]\cite{fellerII},
Steutel and van Harn \cite[\S V.7, Theorems 7.1 and 7.14]{MR2011862},
and Ramachandran and Lau \cite[\S 3.1]{MR1132671} 
for similar treatments.
\begin{proposition}
\label{lmmalpha}
The norming constants $a_n$ associated with a non-degenerate stable distribution are necessarily of the form $a_n = n^{1/\alpha}$ for some $\alpha > 0$.
\end{proposition}

It will be shown by Lemma \ref{lt2} in the next section that necessarily $\alpha \le 2$, yielding Corollary 
\ref{scaling} as a first step towards the ``only if'' part of Theorem \ref{thmstable}.
We will make use of the following elementary lemma.

\begin{lemma}
\label{lmmcf}
If $\phi$ is the characteristic function of a non-degenerate probability distribution, with $|\phi(t)| < 1 $ for some $t$, then
\begin{itemize}
\item[(i)] $| \phi(c t) | \le | \phi(t) |$ for some real $c$ and all real $t$ $\implies |c| \ge 1$.
\item[(ii)] $| \phi(c_1 t)|  = |\phi(c_2 t )|$ for $c_1 >0, c_2 >0$, all real $t$ $\implies c_1 = c_2$.
\end{itemize}

\end{lemma}
\textbf{Proof:} 
It follows from (i) that $|\phi(c^n t ) | \le | \phi(t) |$.
If $|c| < 1$, $c^n \rightarrow 0$ as $n \te \infty$,
and therefore $1 \le |\phi(t)|$. Hence $|c| \ge 1$.

\quad
If the assumption of (ii) holds, $| \phi(t)| = | \phi(c_2 t / c _1)|$. Therefore $c_2/c_1 \ge 1$.
Similarly $c_1/c_2 \le 1$, and so $c_2/c_1 = 1$.
$\square$

\begin{lemma}
\label{lmmanbn}
The constants $a_n$ and $b_n$ associated with a non-degenerate stable distribution of $X$ as in 
{\em \eqref{s1}} or {\em \eqref{s2}} or {\em \eqref{11}}
are such that
\begin{align}
\label{inc}
1 = a_1 < a_2 < \cdots , 
\end{align}
\begin{align}
\label{2}
a_{mn} = a_m a_n, 
\end{align}
\begin{align}
\label{3}
b_m( n - a_n ) = b_n ( m - a_m ).
\end{align}
\begin{align}
\label{nozeros}
\mbox{Moreover, the characteristic function of a stable distribution has no real zeros.}
\end{align}
\end{lemma}

\textbf{Proof:} 
By assumption, the distribution is not degenerate, so $|\phi(t)|$ is not identically equal to $1$.
Observe that
\begin{align*}
\phi(t)^n &= \phi(a_n t ) e^{i b_n t }, \mbox{ all real } t, a_n >0, b_n \mbox{ real}. \\
|\phi(a_{n+1}t)| &= |\phi(t)| \,|\phi(a_n t ) | \\
& \le | \phi(a_n t ) | \\
|\phi(a_{n+1}t/a_n )| & \le |\phi(t)| .
\end{align*}
Thus $a_{n+1} \ge a_n$:
but $a_{n+1} = a_n$ would imply $|\phi(t)| \equiv 1$, and therefore $a_{n+1} > a_n$.
This gives \eqref{inc}.
Next, to see that $\phi$ has no zeros,  observe that
\begin{align*}
|\phi(t)|^2 &= |\phi(a_2 t)|, \\
 |\phi(t/a_2)|^2 &= |\phi(t)|.
\end{align*}
Hence
\begin{align*}
\phi(t) = 0 
\implies
\phi(t/a_2) = 0 
\implies
\phi(t/a_2^n) = 0 
\implies
\phi(0) = 0.
\end{align*}
So $\phi(t) \ne 0$ for real $t$. 

Turning to the proof of \eqref{2} and \eqref{3}, observe that
\begin{align*}
\phi(t)^{mn} &= \phi(a_m t )^n \exp( i n b_m t ) = \phi(a_m a_n t ) \exp i (  n b_m  + a_m b_n) t 
\end{align*}
But also
\begin{align*}
\phi(t)^{mn} &= \phi(a_n a_m t ) \exp i (  m b_n  + a_n b_m) t  = \phi( a_{mn}t ) \exp  i b_{mn} t .
\end{align*}
Equating these expressions gives  for all positive integers $m$ and $n$
\begin{align*}
| \phi( a_{mn}t ) |  = | \phi(a_n a_m t )  |,
\end{align*}
hence $a_{mn} = a_n a_m$ by Lemma \ref{lmmcf} (ii).
Similarly
$$
\exp i (  n b_m  + a_m b_b) t  = \exp i (  m b_n  + a_n b_m) t   ,
$$
hence $ n b_m + a_m b_n = m b_n + a_n b_m$, which is \eqref{3}.
$\square$.

The proof of the Proposition \ref{lmmalpha} stated at the beginning of this section is completed by the following Lemma, which appears
in the same context in the texts of 
Bergstr\"om \cite[\S 8.4]{MR0161363} 
and Ramachandran and Lau \cite[Proposition 1.1.9]{MR1132671}. 
\begin{lemma}
\label{note1}
If a sequence $a_n$ is strictly increasing and satisfies the functional equation $a_{mn} = a_m a_n$, then $a_n = n^k$ for some $k >0$.
\end{lemma}
\textbf{Proof:} 
If $p$ is a positive integer, $a_{m^p} = a_m^p$.
If $1 < m < n$, there is a positive integer $q$ such that
\begin{align*}
n^p & \le m^q  < n^{p+1} \\
a_n^p & \le a_m^q  < a_n^{p+1} \\
p \log a_n & \le q \log a_m  < (p+1) \log a_n \\
(p+1) \log n &> q \log m \ge p \log n \\
\frac{p}{p+1} \frac{ \log a_n }{\log n } &< \frac{ \log a_m} {\log  m } < \frac{ p+1}{p} \frac{ \log a_n } {\log n } 
\end{align*}
When $p \te \infty$ this gives
\begin{align}
\frac{ \log a_n }{\log n } &= \frac{ \log a_m} {\log  m }  = \mbox{ some constant } k ,
\end{align}
\begin{align}
\log a_n  = k \log n,  \quad  a_n = n^k, k >0.
\end{align}
$\square$

\subsection{The case $\alpha \ne 1$}
\label{casene1}

The arguments in this subsection and the next follow an approach developed also by Ramachandran and Lau \cite[Lemma 3.1.4]{MR1132671}. 
Continuing to assume that $\phi$ is a non-degenerate stable characteristic function with associated constants
$a_n$ and $b_n$, we know from Lemma \ref{note1} that $a_n = n^k$ for $k = 1/\alpha >0$. 
Suppose now that $k \ne 1$.
It follows from \eqref{3} that
$$
b_n = \gamma ( n - a_n) = \gamma ( n - n^k )
$$
where $\gamma$ is a real constant. Now from \eqref{1}
\begin{align*}
\phi(t)^n &= \phi( a_n t ) \exp (i b_n t ) \\
n \log \phi(t) &= \log \phi ( n^k t ) + i \gamma ( n - n^k) t \\
n \log ( \phi(t) - i \gamma t ) &=  \log \phi ( n^k t ) - i \gamma  n^k t.
\end{align*}

Put 
\begin{align}
\label{xpsi}
\xpsi(t) = \log \phi(t) - i \gamma t .
\end{align}
Then
\begin{align*}
\xpsi(n^k t )  &= n \, \xpsi(t) \\
\xpsi(n^k t /n^k )  &= n \, \xpsi(t/n^k) \\
\xpsi( t /n^k )  &= \xpsi(t)/n \\
\xpsi \left( \frac{m^k}{n^k} t \right)  &= m \, \xpsi(t/n^k) = \frac{m}{n} \xpsi(t). \\
\end{align*}
Thus for any positive rational $r$, hence also by continuity for all real $r >0$ and $t >0$,
\begin{align*}
\xpsi(r^k t ) &= r \xpsi(t) \\
\xpsi( r t ) &= r ^{1/k}\xpsi(t) = r^{\alpha} \xpsi(t) =  t^\alpha \xpsi(r)
, \mbox{ where } \alpha = 1/k.\\
\end{align*}
Taking $r = 1$ we obtain for $t >0$
\begin{align*}
\xpsi(t) &= t^\alpha \xpsi(1) = - ( c - i c_1) t^\alpha, \mbox{ with } c, c_1 \mbox{ real}
\end{align*}
That is, from the definition \eqref{xpsi}, for $ t >0$
\begin{align*}
\log \phi(t) &= i \gamma t - ( c - i c_1) t^\alpha
\end{align*}
For $t < 0$, by complex conjugation,
$$
\log \phi(t) = \log \overline{ \phi(|t|)} =  \overline{\log \phi(|t|)} =  i \gamma t - ( c + i c_1) |t|^\alpha .
$$
Combine the two cases to obtain
\begin{align}
\label{twosided}
\log \phi(t) = i \gamma t - c |t|^\alpha ( 1 - i c'  \sgn(t)) 
\end{align}
with $c' = c_1/c$.
Here $c >0$ because
$|\phi(t)| = \exp( - c |t|^\alpha) \le 1$, and $c \ne  0$ by the assumption that the distribution is not degenerate.

\begin{lemma}
\label{lt2}
The index $\alpha$ of a stable law is such that $0 < \alpha \le 2$.
\end{lemma}
\textbf{Proof:} 
It is enough to consider the case $\alpha \ne 1$.
Let $U(t)$ be the real part of $\phi(t)$ for real $t$. 
Then \eqref{twosided} gives
$$
1 - U(t) \sim c |t|^\alpha \mbox{ as } t \te 0.
$$
Recall that for the characteristic function of a random variable $X$
$$
\lim_{t \te 0} \frac{ 1 - U(t)}{t^2} = E(X^2), \mbox{ finite or infinite}.
$$
See for instance \cite[(3.8)-(3.9)]{MR2722836}. 
If $\alpha > 2$, this would give $E(X^2) = 0$. Thus $\alpha \le 2$. 
$\square$

See also Feller \cite[p. 171]{fellerII} for another proof that $\alpha \le 2$.

\begin{lemma}
\label{note2lmm}
In the representation {\rm \eqref{twosided}} for $\alpha \ne 1$ the constant $c'$ is subject to the constraint
$$
c' = \beta \tan \pi \alpha/2, \mbox{ where } -1 \le \beta \le 1, \mbox{ if } 0 < \alpha < 2
$$
and $c'= 0$ when $ \alpha = 2$.
\end{lemma}

See Section \ref{constantsec} for the proof.

Combining the above results, we finally have 
\eqref{omeganot1},
$$
\log \phi(t) = i \gamma t  - c |t|^\alpha \left( 1 - i \beta \tan (\pi \alpha /2) \sgn(t) \right)
$$
$$
\mbox{ for } 0 < \alpha \le 2, \alpha \ne 1, -1 \le \beta \le 1, c > 0.
$$
This is the form \eqref{lkstable} for $\alpha \ne 1$, just with $\mu$ replaced by $\gamma$ and $c^\alpha$ replace by $c$.

\subsection{The case $k = 1, a_n = n, \alpha = 1$}
\label{caseeq1}
From \eqref{1}
\begin{align*}
n \log \phi(t)   &= \log \phi(nt) + i b_n t \\
m n \log \phi(t) &= m ( \log \phi(nt) + i b_n t )\\
		 &= \log \phi(mnt) + i b_m n t  + i b_n m t .
\end{align*}
Put $B_n = b_n/n$.
\begin{align*}
m n \log \phi(t) &= \log \phi ( m n t ) + i ( B_m + B_n ) m n t \\
&= \log \phi ( m n t ) + i  B_{mn} t  .
\end{align*}
Thus

\begin{align*}
B_{mn} &= B_{m} + B_{n} \\
n \log \phi(t) &= \log \phi(n t) + i B_n n t \\ 
\frac{1}{n} \log \phi(t) &= \log \phi( t / n ) - i B_n t /n . \\
B_{1/n} &= - B_n
\end{align*}

By similar arguments we can show that
\begin{align*}
B_{m/n} = B_m + B_{1/n} = B_m - B_n .
\end{align*}

This defines $B_r$ for all positive rationals $r$, in such a way that
\begin{align}
r \log \phi(t) &= \log \phi(rt) + i B_r r t   \label{br1}\\
\end{align}
and
\begin{align}
B_{rs } &= B_r + B_s, \quad r >0,  s>0. \label{br2}
\end{align}
By continuity of $\phi$, formula \eqref{br1} serves to define $B_r$ as a continuous function of $r$ for all real $r >0$. The functional equation
\eqref{br2} then extends by continuity from all positive rational $r$ and $s$ to all positive real $r$ and $s$. Let $B(r):= B_r$.
Then for $p,q$ real
$$
B( e^{p+q} ) = B(e ^p) + B(e^q ).
$$
The function $B(e^p)$ is a continuous and linear function of $p$, therefore
\begin{align*}
B(e^p) &= \lambda p, \quad \mbox{ for some real } \lambda  \\
B(r) &= \lambda \log r   \\
r  \log \phi(t) &= \log \phi(r  t ) + i \lambda ( \log r )  r   t .   
\end{align*}

Put $t = 1$ and $\log \phi(1) = -c + i y$. Then for all real $r >0 $
\begin{align*}
\log \phi( r ) &= i y r - c r  - i \lambda r \log r ,\\
 | \phi( r )  | &= e^{-c r } .
\end{align*}
Therefore $c >0$, and we may put $\lambda = c c'$,
\begin{align}
\label{phit1}
\log \phi(t) = i y t - c ( t + i c' t \log t ), \quad t >0 .
\end{align}

To obtain the desired conclusion of \eqref{lkstable}-\eqref{omega1} in this case,
it remains to prove the following Lemma, which will be done in Section \ref{constantsec}.

\begin{lemma}
\label{note4lmm}
In the formula {\rm \eqref{phit1}} for the case $\alpha=1$,
$$c' = 2 \beta /\pi \mbox{ for some }-1 \le \beta \le 1.$$
\end{lemma}

\section{Tail balance for regularly varying distributions}
\label{regvar}
Let $X$ be a real valued random variable governed by a probability measure $P$, with
distribution function $F$. 
To accomodate the Fourier transform framework adopted in this section, 
it is necessary to take $F$ to be the {\em intermediate} distribution function
\begin{equation}
\label{intermed}
F(x) := \frac{1}{2} \left[ P( X < x ) + P(X \le x ) \right]
\end{equation}
as considered by L\'evy \cite[p. 39]{levy54}.
Following the notation of \cite{pitman68}, for $x > 0 $ form the {\em tail sum} and {\em tail difference} functions
\begin{align}
\label{sumdiff}
H(x) := 1 - F(x) + F(-x); \qquad K(x) := 1 - F(x) - F(-x) .
\end{align}
The main result in this section is Theorem \ref{thm2}, which is applied in Section \ref{constantsec} to complete
the characterization of the  stable distributions by showing that the constant $\beta$ appearing in  the previous
Lemmas \ref{note2lmm} and \ref{note4lmm} is cosntrained to lie in the interval $[-1,1]$.
Theorem \ref{thm2}  establishes a general tail balance result $K(x)/H(x) \to \beta \in [-1,1]$ for $F$ with 
characteristic function $\phi(t)$ that is regularly varying with suitable index as $t \to 0$.  This is done using the asymptotic
behaviour of Fourier cosine and sine transforms  of $F$, that is the
the real and imaginary parts of characteristic function
$$
\phi(t) = U(t) + i V(t) = \int_{-\infty}^\infty \cos tx \,dF(x) + i \int_{-\infty}^\infty \sin tx \,dF(x)   .
$$
Equivalently, in terms of the tail functions, there are the formulas obtained by integration by parts \cite[p. 424]{pitman68}
\begin{align}
\label{uvformulas}
\frac{ 1 - U(t)}{t} = \int_0^\infty H(x) \sin \, tx   \, dx ;
\qquad \frac{ V(t)}{t} &= \int_0^\infty K(x) \cos \, tx   \, dx , \quad t \ne 0
\end{align}
and the associated inversion formulas
\begin{align}
H(x) = \frac{2 }{\pi} \int_0^\infty \frac{ 1 - U(t) }{t} \sin \,  t x \, dt , \qquad 
K(x) = \frac{2 }{\pi} \int_0^\infty \frac{ V(t) }{t} \cos \,  t  x \, dt , \quad x >  0  .
\label{invform}
\end{align}
As observed in \cite[\S 8.1.4, p. 336]{bgt}, these inversion formulas follow easily, by linear combination of components using \eqref{sumdiff},
from the inversion formula of Gil-Paez \cite{MR0045992}
\begin{align}
\label{gilpaez}
F(x) 
= \frac{1}{2} -  \frac{1} { \pi } \int_0^ \infty   \frac{ \Im [ e^{-itx} \phi(t) ] } {t} \, dt
= \frac{1}{2} +  \frac{1} {2 \pi } \int_0^ \infty  \frac{ e^{itx} \phi(-t) - e^{-itx} \phi(t)  }{i t} \, dt
\end{align}
in which the choice of the intermediate version \eqref{intermed} of the distribution function $F$ is essential if $F$ has discontinuities.
For degenerate $F$ the formula \eqref{gilpaez} reduces to the Dirichlet integral \eqref{dirichlet}, and 
the general case of \eqref{gilpaez} follows by a Fubini argument.  Note that the integrals involved in these inversion formulas 
\eqref{gilpaez} and \eqref{uvformulas} are improper, 
obtained as limits of $\int_\epsilon^T$ as $\epsilon \downarrow0$ and $T \uparrow \infty$.
See \cite{MR0120676} for further discussion. 
The better known inversion formula of 
L\'evy  \cite[p. 38]{levy54}
\begin{align}
\label{levy}
F(y) - F(x)
= \frac{1}{2 \pi} 
\lim_{T \uparrow \infty} \int_{-T}^T  \frac{ ( e^{- itx } - e^{-ity} ) \phi(t)  }{i t} \, dt
\end{align}
is the difference of two evaluations of \eqref{gilpaez}, which can also be derived quite easily from \eqref{levy} \cite[p. 39]{levy54}.

The relations between asymptotic behaviour of $H$ and $K$ at $\infty$ and of $U$ and $V$ at $0$ involve the basic cosine and sine integral evaluations
$C(k)$ and $S(k)$ of Corollary \ref{sincos}.  That is
\begin{equation}
\label{sk}
S(k) := \Gamma(1 - k) \cos  \pitwo k = \int_0^\infty \frac{ \sin x }{x^k} dx \mbox { for } 0 < k < 2
\end{equation}
where $S(1) = \pi/2$ is the Dirichlet integral, and

\begin{align}
C(k) := \Gamma(1 - k) \sin  \pitwo k &= \int_0^\infty \frac{ \cos x }{x^k} dx \mbox { for } 0 < k < 1 \label{ck}\\
&= \int_0^\infty \frac{ \cos x - 1 }{x^k} dx \mbox { for } 1 < k < 3.
\end{align}
The following two lemmas were established in \cite{pitman68}. 
See also \cite[\S 8.1.4, p. 336]{bgt} and \cite{MR0361644} for various generalizations.

\begin{lemma}
\label{lemmhu}
{\rm \cite{pitman68}}
If $H(x)$ is regularly varying at $\infty$ with index $-k$, $0 < k < 2$,
i.e.
$$
\mbox{ for all } \lambda > 0, \quad H(\lambda x) / H(x) \te \lambda ^{-k} \mbox{ as } x \te \infty, 
$$
then
\begin{align}
\label{4}
1 - U(t) \sim S(k) H(1/t) \mbox{ as } t \downarrow 0 .
\end{align}
Conversely, if $1-U(t)$ is of index $k$ at $0$, $0 \le k < 2$, then
\begin{align}
\label{5}
H(x) \sim \frac{ 1 - U(1/x) }{S(k) } \mbox{ as } x \te \infty .
\end{align}
\end{lemma}
Here the notation $\sim$ is used for {\em asymptotic equivalence}, meaning that the ratio of two expressions tends to $1$ in the
specified limit regime.
\begin{lemma}
\label{lemmkv}
{\rm \cite{pitman68}}
If $K(x)$ is regularly varying with index $-k$ at $\infty$, 
and if also $K(x)$ is monotonic when $x$ is large, then
\begin{align}
\label{6}
V(t) \sim C(k) K(1/t) \mbox{ as } t \downarrow 0 
\qquad \mbox{ if }   0 < k < 1, 
\end{align}
whereas if $1 < k < 3$, then the distribution has a finite mean $\mu$, and
\begin{align}
\label{7}
V(t) - \mu t \sim C(k) K(1/t) \mbox{ as } t \downarrow 0 \qquad \mbox{ if }   1 < k < 3.
\end{align}
\end{lemma}
See \cite[p. 440, Theorem 8]{pitman68}. (There is a misprint on p. 441. The second line should read $2 n - 1 < m < 2 n + 1$.)

For the proof of Theorem \ref{thm2} below we need the following unsurprising theorem, which is not given in \cite{pitman68}, nor,
as far as we know, anywhere else in the literature.
\begin{theorem}
\label{thm1}
If $V(t)$ is regularly varying at $0$ with index $k$, $0 < k < 2$, $k \ne 1$, and if also $K(x)$ is monotonic when $x$ is great, then
\begin{align}
\label{thm1concl}
K(x) \sim V(1/x) / C(k) \mbox{  as } x \te \infty.
\end{align}
\end{theorem}
This is established in much the same way as the converse part of Lemma \ref{lemmhu} for $1-U(t)$ and $H(x)$. 
See \cite[p.432]{pitman68}.  The details of this argument are provided  in Section \ref{regvarproof}.

\begin{theorem}
\label{thm2}
If $1 - U(t)$ is regularly varying at $0$ with index $k$, $0 < k < 2$, $k \ne 1$, and
$$
V(t) \sim c' ( 1 - U(t) ) \mbox{ as } t \downarrow 0,
$$
then
\begin{equation}
\label{cprime}
c' = \beta \tan \pi  k  /2, \mbox{ where } \beta = \lim_{x \to \infty} \frac{ K(x) }{ H(x) }  \in [ -1,1] .
\end{equation}
\end{theorem}

\textbf{Proof:} 
By application of Lemma \ref{lemmhu}, the assumption on $1- U(t)$ gives
\begin{align}
H(x) &\sim \frac{ 1 - U(1/x) }{ S(k)} \mbox{ as } x \te \infty. \\
1 - U(t) &\sim S(k) H(1/t) \mbox{ as } t \downarrow 0. \\
V(t) &\sim c' ( 1 - U(t) ) \sim c' S(k) H(1/t) \mbox{ as } t \downarrow 0 \label{vdown}. 
\end{align}
In the following analysis we will apply the preceding results with both $|X|$ and 
$X_{+}:= \max(X,0)$ in place of $X$. However, the notations $F$, $K$ and $H$ should always 
be understood relative to the original distribution of $X$.
Let
\begin{equation}
\label{wdef}
W(t) := t \, \int_0^\infty H(x) \cos \, tx \, dx.
\end{equation}
Observe from \eqref{sumdiff} and \eqref{uvformulas} that $W(t)$ is the imaginary part of the characteristic function of $|X|$. 
for $X$ with characteristic function $U(t) + i V(t)$. So by 
Lemma \ref{lemmkv} applied to $|X|$ instead of $X$ 
\begin{equation}
\label{wasym}
W(t) \sim C(k) H(1/t) \mbox{ as } t \downarrow 0 \mbox{ if } 0 < k < 1 
\end{equation}
and 
\begin{equation}
\label{wasym2}
W(t) - E|X| \,t \sim C(k) H(1/t) \mbox{ as } t \downarrow 0 \mbox{ if } 1 < k <  2.
\end{equation}
For $x > 0$ consider the right hand tail probability function
\begin{align}
\label{rdef}
R(x) := 1 - F(x) = \hf ( H(x) + K(x) ), \mbox{ which is monotonic}.
\end{align}
Then, by \eqref{uvformulas} and \eqref{wdef},
the function 
\begin{align*}
 t \, \int_0^\infty R(x) \cos \, tx \, dx = \hf (  W(t) + V(t) )
\end{align*}
is the imaginary part of the characteristic function of $X_{+}:= \max(X,0)$.
Supposing first that $0 < k < 1$, we see from \eqref{vdown} and \eqref{wasym} that
\begin{equation}
\label{wvasym}
\hf (  W(t) + V(t) ) \sim \hf ( C(k) + c' S(k) ) H(1/t) \mbox{ as } t \downarrow 0.
\end{equation}
By Theorem \ref{thm1} applied to $X_{+}$ instead of $X$,
\begin{equation}
\label{rasym}
R(x) \sim \frac{ ( C(k) + c' S(k) ) H(x) }{ 2 C(k) }  = \hf \left( 1 + \frac{ c' }{\tan \, k \pi/2 } \right) H(x) \mbox{ as } x \te \infty .
\end{equation}
But $0 \le R(x)/H(x) \le 1$ by \eqref{sumdiff} and \eqref{rdef}, so the conclusion \eqref{cprime} follows for $0 < k < 1$.
Also if $1 < k < 2$, the same relation \eqref{rasym} holds, leading to the same conclusion.
To see this, observe that for $1 < k < 2$ the assumptions of the theorem imply that $E|X| < \infty$ and $E(X) = 0$, so $E(X_+) = \hf E|X|$.
It follows from \eqref{wasym2} that \eqref{wvasym} still holds provided a term 
$E(X_+) t$ is subtracted from the  left side, which implies the same asymptotics for the imaginary part 
of the characteristic function of the centered random variable $X_+ - E(X_+)$. Theorem \ref{thm1} applied to this centered variable 
then gives \eqref{rasym}, at first with $R(x + E(X_+) )$ instead of $R(x)$, but then also without the shift, by application of
the regular variation of $H(x)$.
$\square$

\section{Proof of Theorem \ref{thm1}}
\label{regvarproof}

We start from the inversion formula \eqref{invform} 
\begin{align}
K(u) = \frac{2 }{\pi} \int_0^\infty \frac{ V(t) }{t} \cos \, u t \, dt , \quad u \ge 0 . \label{invKu} 
\end{align}
We require the following results, obtained by integrating under the integral sign with respect to $u$.
\begin{lemma}
If $V(t)$ is of one sign in the right hand neighbourhood of $0$, then
\begin{align}
K_1(x)&:= \int_0^x K(u) du = \frac{2 }{\pi} \int_0^\infty \frac{ V(t) \sin xt } {t^2} dt , \label{k1}\\
K_2(x)&:= \int_0^x K_1(u) du = \frac{2 }{\pi} \int_0^\infty \frac{ V(t) (1 - \cos \, x t )} {t^3} dt . \label{k2} 
\end{align}
\end{lemma}
\textbf{Proof:}
The inversion formula \eqref{invKu} gives
$$
K_1(x) = \frac{2 }{\pi} \int_0^x du \int_0^\infty \frac{ V(t) }{t} \cos \, ut \, dt .
$$
Because $V(t)$ is of one sign in the right hand neighbourhood of $0$, and $\cos xt \te 1$ as $t \te 0$,
$V(t)/t$ is integrable over $(0,T)$ for each $0 < T < \infty$. So
\begin{align*}
\int_0^x du \int_0^T \frac{ V(t)}{t} \cos \, ut  \, dt &= \int_0^T dt \int_0^x \frac{ V(t)}{t} \cos \, ut  \, du  \\
					 	       &= \int_0^T \frac{ V(t)}{t^2} \sin\, xt  \, dt   .
\end{align*}
If we can show that $\int_0^T \frac{ V(t)}{t} \cos \, ut  \, dt $ is bounded for $0 < T < \infty$, $0 < u < x $, this will prove \eqref{k1}.
Now $\int_0^1 \frac{ V(t)}{t} \cos \, ut  \, dt $ is bounded for all $u >0$, and 
$$
\int_1^T  \frac{ V(t)}{t} \cos \, ut  \, dt = \int_1^T  dt \int_{-\infty}^\infty \frac{ \sin vt \, \cos u t }{t} d F(v)
$$
The modulus of the integrant is at most $1$. Therefore
$$
\int_1^T  \frac{ V(t)}{t} \cos \, ut  \, dt = \int_{-\infty}^\infty d F(v) \int_1^T  \frac{ \sin vt \, \cos u t }{t} dt 
$$
where
\begin{align}
\int_1^T  \frac{ \sin vt \, \cos u t }{t} dt &= \int_1^T  \frac{ \sin (v + u ) t  + \sin (v-u) t }{2t} dt  \\
&= \int_{v+u} ^{(v+u) T}   \frac{ \sin  t  }{2t} dt  + \int_{v-u} ^{(v-u) T}   \frac{ \sin  t  }{2t} dt,   \\
\end{align}
which is bounded for all $T,v,u$. The integral $\int_1^T \frac{V(t)}{t} \cos \, ut \,  dt $ will be bounded, hence so will be the $\int_0^T$,
and \eqref{k1} follows. The result \eqref{k2} is then obvious, since we start with an integrand which is absolutely integrable with respect to $t$ over $0 < t < \infty$.
$\square$

\textbf{Proof of Theorem \ref{thm1}:}
It is assumed that $K(x)$ is ulimately monotonic, and $V(t)$ is regularly varying with index $k$, $0 < k < 2$.
So $V(t)$ will be of one sign in some right hand neighbourhood of $0$, and
$$
\frac{K_2(x)}{ x^2 V(1/x)} = \frac{2}{\pi} \int_0^\infty \frac{ V(t/x)}{V(1/x)} \frac{ 1 - \cos t } {t ^3 } \, dt .
$$
When $x \te \infty$, $V(t/x)/V(1/x) \te t^k$, and the integrand converges to $t^{3-k} ( 1 - \cos \, t )$.
Since $V(1/x)$ is of index $-k$ at $\infty$, corresponding to any $h >0$, there exist $A>0, B>0$ such that when $x > B$,
\begin{align*}
\left| \frac{V(t/x)}{V(1/x)} \right| &< A t^{k+h} \mbox{ when } t > 1 \\
& <  A t^{k-h} \mbox{ when } t < 1 .
\end{align*}
See \cite[p. 426, Lemma 2]{pitman68}.
Take $h = 1 - k/2$, so that $k +h = 1 + k/2$.
When $x > 1$, the modulus of the integrand will be at most $A ( 1 - \cos \, t )/ t^{2 - k/2 }$, which is integrable over $(1,\infty)$. 
Therefore, when $x \te \infty$,
$$
\int_1^\infty \frac{ V(t/x)}{V(1/x)} \frac{ 1 - \cos t } {t ^3 } \, dt  \te \int_1^\infty  \frac{ 1 - \cos t } {t ^{3-k} } \, dt . 
$$
Similarly, taking $h = k/2$, we can show that
$$
\int_0^1 \frac{ V(t/x)}{V(1/x)} \frac{ 1 - \cos t } {t ^3 } \, dt  \te \int_0^1 \frac{ 1 - \cos t } {t ^{3-k} } \, dt . 
$$
Thus as $x \te \infty$
\begin{align}
\label{note3formula}
\frac{K_2(x)}{ x^2 V(1/x)}  \te  \frac{2 }{\pi} \int_0^\infty \frac{ 1 - \cos t } {t ^{3-k} } \, dt  = \frac{2 }{\pi} C(3-k) .
\end{align}
Therefore
$$
K_2(x) \sim \frac{2 }{\pi} C(3-k) x^2 V(1/x),
$$
and so is regularly varying at $\infty$ with index $2 -k$. Also, $K(x) \te 0$ as $x \te \infty$ and is ultimately monotonic, so is ultimately of one sign.
So $K_1(x):= \int_0^x K(u) \, du$ is ultimately monotonic and $K_2(x):= \int_0^x K_1(u) du$ is regularly varying with index $2-k$. It follows \cite[p. 427, Lemma 3]{pitman68}
that $K_1(x)$ is regularly varying at $\infty$ with index $1-k$. 

If $0 < k < 1$, a second application of the same Lemma shows that $K(x)$ is regulary varying at $\infty$ with index $-k$.
The required conclusion \eqref{thm1concl} then follows from \eqref{6} above.

If $1 < k < 2$, then $V(t)$ will be of index $>1$, and the distribution will have a finite mean value $K_1(\infty) = \int_0^\infty K(u) du$. This mean value must be zero.
Otherwise, $V(t)$ would have a term $\mu t$, and be of index $1$.
So $\int_0^\infty K(x) dx = 0$ and
$$
K_1(x):= \int_0^x K(u) du = - \int_x^\infty K(u) du \mbox{ regularly varying at $\infty$ with index } 1 - k < 0.
$$
The Lemma quoted above will not apply. However, the Lemma below is proved in the same way as the Lemma quoted above, and will imply that $K(x)$ is of index $-k$ at $\infty$.
The conclusion \eqref{thm1concl} will then follow from \eqref{7}.
$\square$

\begin{lemma}
If $G(x)$ is monotonic when $x$ is great, and $\int_x^\infty G(u) du $ is regularly varying with index $k' <0$ at $\infty$, then $G(x)$ is regularly varying with index $k'-1$.
\end{lemma}

When $k=1$, the function $K_1(x)$ will be of index $0$, and neither Lemma will apply. 
In fact, Theorem \ref{thm1} is not true for $k=1$.

\section{Restrictions on the constants}
\label{constantsec}

\textbf{Proof of Lemma \ref{note2lmm}.}

Starting from the representation \eqref{twosided}, for a stable distribution of $X$ with characteristic exponent $\alpha$, with $0 < \alpha < 2$ and $\alpha \ne 1$,
$$
\log \phi(t) = i \gamma  t - c ( 1 - i c') t^\alpha, \quad t > 0.
$$
By working with $X - \gamma$ instead of $X$, we can make $\gamma = 0$. Then
\begin{align*}
\phi(t) = \exp( - c t^\alpha ) ( \cos \, c c' t^\alpha + i \sin \, cc't^\alpha ), \quad t > 0
\end{align*}
which makes
\begin{equation}
\label{uvalpha}
1- U(t) \sim c t^\alpha  \mbox{ and } V(t) \sim c c' t^\alpha \mbox{ as } t \downarrow 0.
\end{equation}
As $V$ is regularly varying with index $\alpha$, Theorem \ref{thm2} applies for $0 < \alpha < 1$
with $\alpha \ne 1$.

The preceding argument does not cover the case $\alpha = 2$. In that case
\begin{align*}
\log \phi(t) &= - c ( 1 - i c') t^2, \quad t >0, \\
\phi(t) &= \exp( - c ( 1 - i c') t^2 ),
\end{align*}
which gives 
\eqref{uvalpha} with $\alpha = 2$.
It follows that the distribution has a finite variance $2 c$, and zero mean. If $X_1, \ldots, X_n$ are
independent random variables with this distribution,
$(X_1 + \cdots + X_n)/\sqrt{ 2 c n }$ has a limit distribution which is standard normal:
as $n \te \infty$
\begin{align*}
\phi( t / \sqrt{ 2 c n } )^n &\te \exp( - \hf t^2 ) \\
n \log \phi( t / \sqrt{ 2 c n } ) &\te  - \hf t^2  \\
\frac{ - n c ( 1 - i c' ) t^2 }{ 2  c n  } &= - \hf t^2 ( 1 - i c' ) \te - \hf t^2.
\end{align*}
Therefore, $c'= 0$. This is consistent with $c' = \beta \tan \alpha/2 = \beta \tan \pi = 0$.
$\square$

\textbf{Proof of Lemma \ref{note4lmm}.}

This is the case $k=1$. From \eqref{phit1} we have
\begin{align*}
\log \phi(t) &= i y t - c (t + i c' t \log t ), \quad t >0. \\
\phi(t) &= \exp( i y t - c t ) ( \cos ( c c' t \log t ) - i \sin ( c c' t \log t ) ), \quad t > 0,
\end{align*}
which makes
$$
1 - U(t) \sim ct \mbox{ and }
V(t) \sim - cc' t \log t \mbox{ as } t \downarrow 0.
$$
As shown in \eqref{note3formula}
\begin{align*}
\frac{K_2(x)}{ x^2 V(1/x)}  &= \frac{2}{\pi} \int_0^\infty \frac{ V(t/x)}{V(1/x)} \frac{ 1 - \cos t } {t ^3 } \, dt . \\
			    &\te  \frac{2 }{\pi} \int_0^\infty \frac{ 1 - \cos t } {t ^{2} } \, dt  = 1 \mbox{ as } x \te \infty
\end{align*}
using \eqref{cm} for $\mm = 2$. Hence 
\begin{align}
\label{k2asym}
K_2(x) \sim x^2 V(1/x) \sim c c' x \log x \mbox{ as } x \te \infty. 
\end{align}
On the other hand, Lemma \ref{lemmhu} gives
\begin{align*}
H(x) \sim  \frac{ 1 - U(1/x) } {S(1)} \sim \frac{ 2 c }{\pi x } \mbox{ as } x \te \infty
\end{align*}
hence
\begin{align*}
H_1(x) := \int_0^x H(u) du \sim \frac{ 2 c }{\pi } \, \log x   \mbox{ as } x \te \infty 
\end{align*}
and
\begin{align}
\label{h2asym}
H_2(x) := \int_0^x H_1(u) du \sim \frac{ 2 c }{\pi } \, x \log x   \mbox{ as } x \te \infty  .
\end{align}
Moreover, $|K(x)| \le H(x)$ by \eqref{sumdiff} and the triangle inequality. Therefore, $|K_2(x)| \le H_2(x)$. So
\eqref{k2asym} and \eqref{h2asym} imply
$ | c c' | \le \frac{2 c }{\pi}$, that is $c' = 2 \beta/\pi$ where $ -1 \le \beta \le 1$.
$\square$

\section*{Acknowledgements}
Thanks to the referee for a careful reading of the paper and for a number of constructive suggestions for improvement.
Thanks also to Simon Ruijsenaars for catching a number of slips in Section 3.

\def\polhk#1{\setbox0=\hbox{#1}{\ooalign{\hidewidth
  \lower1.5ex\hbox{`}\hidewidth\crcr\unhbox0}}}
  \def\polhk#1{\setbox0=\hbox{#1}{\ooalign{\hidewidth
  \lower1.5ex\hbox{`}\hidewidth\crcr\unhbox0}}}
  \def\polhk#1{\setbox0=\hbox{#1}{\ooalign{\hidewidth
  \lower1.5ex\hbox{`}\hidewidth\crcr\unhbox0}}}
  \def\polhk#1{\setbox0=\hbox{#1}{\ooalign{\hidewidth
  \lower1.5ex\hbox{`}\hidewidth\crcr\unhbox0}}}
  \def\polhk#1{\setbox0=\hbox{#1}{\ooalign{\hidewidth
  \lower1.5ex\hbox{`}\hidewidth\crcr\unhbox0}}}
  \def\polhk#1{\setbox0=\hbox{#1}{\ooalign{\hidewidth
  \lower1.5ex\hbox{`}\hidewidth\crcr\unhbox0}}}
  \def\polhk#1{\setbox0=\hbox{#1}{\ooalign{\hidewidth
  \lower1.5ex\hbox{`}\hidewidth\crcr\unhbox0}}} \def\cprime{$'$}
  \def\polhk#1{\setbox0=\hbox{#1}{\ooalign{\hidewidth
  \lower1.5ex\hbox{`}\hidewidth\crcr\unhbox0}}}
  \def\polhk#1{\setbox0=\hbox{#1}{\ooalign{\hidewidth
  \lower1.5ex\hbox{`}\hidewidth\crcr\unhbox0}}}
  \def\cftil#1{\ifmmode\setbox7\hbox{$\accent"5E#1$}\else
  \setbox7\hbox{\accent"5E#1}\penalty 10000\relax\fi\raise 1\ht7
  \hbox{\lower1.15ex\hbox to 1\wd7{\hss\accent"7E\hss}}\penalty 10000
  \hskip-1\wd7\penalty 10000\box7}
  \def\polhk#1{\setbox0=\hbox{#1}{\ooalign{\hidewidth
  \lower1.5ex\hbox{`}\hidewidth\crcr\unhbox0}}}
  \def\polhk#1{\setbox0=\hbox{#1}{\ooalign{\hidewidth
  \lower1.5ex\hbox{`}\hidewidth\crcr\unhbox0}}}
  \def\cftil#1{\ifmmode\setbox7\hbox{$\accent"5E#1$}\else
  \setbox7\hbox{\accent"5E#1}\penalty 10000\relax\fi\raise 1\ht7
  \hbox{\lower1.15ex\hbox to 1\wd7{\hss\accent"7E\hss}}\penalty 10000
  \hskip-1\wd7\penalty 10000\box7}
  \def\polhk#1{\setbox0=\hbox{#1}{\ooalign{\hidewidth
  \lower1.5ex\hbox{`}\hidewidth\crcr\unhbox0}}} \def\cprime{$'$}
  \def\polhk#1{\setbox0=\hbox{#1}{\ooalign{\hidewidth
  \lower1.5ex\hbox{`}\hidewidth\crcr\unhbox0}}} \def\cprime{$'$}
  \def\cprime{$'$} \def\cprime{$'$} \def\cprime{$'$} \def\cprime{$'$}
  \def\cprime{$'$} \def\cprime{$'$}
  \def\polhk#1{\setbox0=\hbox{#1}{\ooalign{\hidewidth
  \lower1.5ex\hbox{`}\hidewidth\crcr\unhbox0}}}
  \def\polhk#1{\setbox0=\hbox{#1}{\ooalign{\hidewidth
  \lower1.5ex\hbox{`}\hidewidth\crcr\unhbox0}}}
  \def\polhk#1{\setbox0=\hbox{#1}{\ooalign{\hidewidth
  \lower1.5ex\hbox{`}\hidewidth\crcr\unhbox0}}}
  \def\polhk#1{\setbox0=\hbox{#1}{\ooalign{\hidewidth
  \lower1.5ex\hbox{`}\hidewidth\crcr\unhbox0}}}
  \def\polhk#1{\setbox0=\hbox{#1}{\ooalign{\hidewidth
  \lower1.5ex\hbox{`}\hidewidth\crcr\unhbox0}}}
  \def\Dbar{\leavevmode\lower.6ex\hbox to 0pt{\hskip-.23ex \accent"16\hss}D}
  \def\polhk#1{\setbox0=\hbox{#1}{\ooalign{\hidewidth
  \lower1.5ex\hbox{`}\hidewidth\crcr\unhbox0}}}
  \def\cfgrv#1{\ifmmode\setbox7\hbox{$\accent"5E#1$}\else
  \setbox7\hbox{\accent"5E#1}\penalty 10000\relax\fi\raise 1\ht7
  \hbox{\lower1.05ex\hbox to 1\wd7{\hss\accent"12\hss}}\penalty 10000
  \hskip-1\wd7\penalty 10000\box7}
  \def\polhk#1{\setbox0=\hbox{#1}{\ooalign{\hidewidth
  \lower1.5ex\hbox{`}\hidewidth\crcr\unhbox0}}}
  \def\polhk#1{\setbox0=\hbox{#1}{\ooalign{\hidewidth
  \lower1.5ex\hbox{`}\hidewidth\crcr\unhbox0}}}
  \def\polhk#1{\setbox0=\hbox{#1}{\ooalign{\hidewidth
  \lower1.5ex\hbox{`}\hidewidth\crcr\unhbox0}}} \def\cprime{$'$}
  \def\polhk#1{\setbox0=\hbox{#1}{\ooalign{\hidewidth
  \lower1.5ex\hbox{`}\hidewidth\crcr\unhbox0}}}
  \def\polhk#1{\setbox0=\hbox{#1}{\ooalign{\hidewidth
  \lower1.5ex\hbox{`}\hidewidth\crcr\unhbox0}}}
  \def\polhk#1{\setbox0=\hbox{#1}{\ooalign{\hidewidth
  \lower1.5ex\hbox{`}\hidewidth\crcr\unhbox0}}}
  \def\polhk#1{\setbox0=\hbox{#1}{\ooalign{\hidewidth
  \lower1.5ex\hbox{`}\hidewidth\crcr\unhbox0}}}

\end{document}